Su dongcai

# Stability analysis of delay differential equations via Semidefinite programming

Dongcai Su

## Abstract

This paper study the problem of how to choose parameter $\theta$ so that the delay differential equation (DDE): $\frac{d}{dt}x(t) = A_0(\theta)x(t) + \sum_{i=1}^{K} A_i(\theta)x(t-\tau_i)$ is stable? where $x(t)$ denotes the $n$ dimensional vector, $A_i(\theta), 0 \leq i \leq K$ are $n \times n$ matrices whose entries are depend linearly on the parameter $\theta$, $\tau_i > 0, 1 \leq i \leq K$ are positive numbers indicating the time delays. After discretizing the DDE above, we show that the problem can be casted equivalently into a semi-definite programming (SDP) see (3.2), which can be solved effectively through some popular method, e.g., the interior point method [1].

## 0. Introduction

To analyze a linear delay differential equations (DDE) in general is difficult, consider a typical form of DDE, see equation (0.1),

$$\frac{d}{dt}x(t) = A_0(\theta)x(t) + \sum_{i=1}^{K} A_i(\theta)x(t-\tau_i) \quad (0.1)$$

Here in (0.1), $x(t) = \left[x_1(t), \cdots, x_n(t)\right]^T$ is a n-dimensional vector, $\frac{dx(t)}{dt} = \left[\frac{dx_1(t)}{dt}, \cdots, \frac{dx_n(t)}{dt}\right]^T$, $A_i(\theta), 0 \leq i \leq K$ are $n \times n$ matrices whose entries depend linearly on the parameter $\theta$, $\tau_i > 0, 1 \leq i \leq K$ are positive numbers indicating the time delays.

According to [2], the character equation of (0.1) is

$$\det\left(-\lambda I + A_0(\theta) + \sum_{i=1}^{K} A_i(\theta)\exp(-\tau_i \lambda)\right) \quad (0.2)$$

The stability of (0.1) is therefore depends on the solutions $\lambda$ of (0.2), however, since (0.2) contains the exponential term, which in general can have infinite many solutions and therefore difficult to solve analytically, even if it's solved in a closed form, the solutions might be too complicated to draw any useful conclusions on its behavior [3].

However, understanding and analyzing the stability of DDE like (0.1) is important in practice, for example, in engineering, the time delay factors in (0.1) reflect human reaction time or device responding time, and parameter $\theta$ is associate to some control parameters which we can manipulate, therefore a crucial problem in this case is to turn the parameter $\theta$ in real time, so that systems like (0.1) is stable.

Motivated by this, this paper is dedicated to answer the question: "How to choose $\theta$, so that (0.1) is stable, or declare no such $\theta$ exists." To this end, first we discretize (0.1) into (2.1), and then after arguing with Schur complement, we show that the stability of (2.1) is equivalent to the existence of a solution of SDP (3.2) with $t<1$.

The remaining part of this paper is organized as following, in section 1 we state the problem studied in this paper, which is the stability of DDE (0.1), in section 2 we propose a discretized version of (0.1), which is (2.1) and argue that stability conditions for (0.1) and (2.1) are equivalent, finally, in section 3, we analyze the stability condition of (2.1) by a SDP approach.

# 1. Problem statement

The problem studied in this paper is "How should we choose the parameter $\theta$, so that the delayed differential equations system (DDE) in (0.1) is stable? " Here, "DDE (0.1) is stable" we mean $\lim_{t\to\infty}\|x(t)\|_2 = 0$, here $\|\cdot\|_2$ denotes the Euclidean norm.

# 2. Interpret (0.1) in a discretized version

To analyze the stability of (0.1), we approximate it with a discretized version, before we present the discretized version of (0.1), we introduce some notations as explained below, let
$u(t) = \left[x_1(t), x_1(t-\delta t), \cdots, x_1(t-N\delta t), \cdots, x_n(t), x_n(t-\delta t), \cdots, x_n(t-N\delta t)\right]^T$ which denotes a $nN$ dimensional vector, here $\delta t > 0$ is a small positive number denotes the time step, $N\delta t = \max\{\tau_1, \cdots, \tau_K\}$ denotes the maximum time delay[1]. Then after applying a discretized version of derivative based on the technique introduced in appendix A, we can interpret (0.1) as,

$$\frac{1}{\delta t}\left(Du(t) + D^{(r)}u(t-N\delta t)\right) = B(\theta)u(t) + A(\theta)u(t-N\delta t) \qquad (2.1)$$

Where in the right hand side of (2.1), $B(\theta)$, $A(\theta)$ are matrices with proper size whose

---
[1] For simplicity, we assume $\tau_i$, $1 \leq i \leq K$ in (0.1) are positive rational numbers, so that entries of u(t) contains all entries of x(t − $\tau_i$), $1 \leq i \leq K$ if we choose $\delta t$ sufficiently small.

entries are dependent linearly on parameter $\theta$. The left hand side of (2.1) is nothing but discretized version of $\frac{d}{dt}u(t)$, where $D = \begin{bmatrix} D_1 & & \\ & \ddots & \\ & & D_n \end{bmatrix}$, $D^{(r)} = \begin{bmatrix} D_1^{(r)} & & \\ & \ddots & \\ & & D_n^{(r)} \end{bmatrix}$, here

$D_i = \begin{bmatrix} w_1, \cdots, w_m & & \\ & \ddots & \\ & & w_1 \end{bmatrix}$, $1 \leq i \leq n$ are identical, with,

$$D_i(j,k) = \begin{cases} w_{k-j+1}, & \text{if } k \geq j, 1 \leq k-j+1 \leq m \\ 0, & \text{otherwise} \end{cases}, 1 \leq j \leq N, 1 \leq k \leq N \quad (2.2)$$

$D_i^{(r)} = \begin{bmatrix} 0 & & 0 \\ w_m & & \\ \vdots & \ddots & 0 \\ w_2 & \cdots & w_m \end{bmatrix}$, $1 \leq i \leq n$ are identical, with,

$$D_i^{(r)}(j+n-m+1, k) = \begin{cases} w_{m-(j-k)}, & \text{if } k \leq j, 1 \leq j \leq m-1, 1 \leq k \leq m-1 \\ 0, & \text{otherwise} \end{cases}, \quad (2.3)$$

for $1 \leq j \leq N, 1 \leq k \leq N$

For convenience, we call the matrix defined in (2.2) as the differentiate matrix, and the matrix defined in (2.3) as differentiate remaining matrix.

Reformulating (2.1) as,

$$u(t) = (D - \delta t B(\theta))^{-1} (-D^{(r)} + \delta t A(\theta)) u(t - N\delta t) \quad (2.4)$$

we can see that (2.1) is stable if and only if,

$$\left\| (D - \delta t B(\theta))^{-1} (-D^{(r)} + \delta t A(\theta)) \right\|_{2 \to 2} < 1 \quad (2.5)$$

Here in (2.5), $\|\cdot\|_{2 \to 2}$ denotes the operator norm, which is the maximal eigen-value (measured in absolute value) of a matrix. (2.5) is equivalent to,

$$(-D^{(r)} + \delta t A(\theta))^T (D - \delta t B(\theta))^{-T} (D - \delta t B(\theta))^{-1} (-D^{(r)} + \delta t A(\theta)) \prec I \quad (2.6)$$

By Shcur complement, (2.6) is rewritten as,

$$\begin{bmatrix} (D - \delta t B(\theta))(D - \delta t B(\theta))^T, & -D^{(r)} + \delta t A(\theta) \\ (-D^{(r)} + \delta t A(\theta))^T, & I \end{bmatrix} \succ 0 \quad (2.7)$$

Applying Shcur complement of $I$ in (2.7) shows that (2.7) is equivalent to,

$$(D - \delta t B(\theta))(D - \delta t B(\theta))^T - (-D^{(r)} + \delta t A(\theta))^T (-D^{(r)} + \delta t A(\theta)) \succ 0 \quad (2.8)$$

For notational simplicity, we rewrite (2.8) as,

$$DD^T - D^{(r)T}D^{(r)} + \delta t E(\theta) + (\delta t)^2 F(\theta) \succ 0 \qquad (2.9)$$

Where $E(\theta)$ is a matrix of proper size whose entries depend linearly to parameter $\theta$, and $F(\theta)$ is a matrix of proper size whose entries may be non-linear function of parameter $\theta$. $DD^T - D^{(r)T}D^{(r)}$ is the constant matrix term in (2.9) which is introduced by the differential operator on $u(t)$.

**Remarks:**

1. Since according to Appendix A, the derivative of $x_i(t), 1 \leq i \leq n$ can be approximated with a linear combination of its historical data $x_i(t-\delta t), \ldots, x_i(t-k\delta t)$ with the approximation error asymptotically bounded by $O(\delta t)^{k+1}$. Notice that in the discretized system (2.1), constants $w_1, \ldots, w_m$ in matrix $D$, $D^{(r)}$ are nothing but the linear combination coefficients as suggested in Appendix A, therefore after choosing these constants appropriately, (2.1) can faithfully reflect the dynamic of the original DDE system (0.1), in other words, (0.1) is stable if and only if (2.1) is stable provided that the time step $\delta t$ is sufficiently small. Therefore in the remaining part of this article, we dedicate to analyze the stability of (2.1).

2. Here our discretization technique which approximate (0.1) with (2.1) is also refer to as the method of time steps [4, 5] in literatures, which is frequently adopted to obtain the numerical solution of DDE like (0.1) see [6] for instance, and recently, it's applied also in analyzing the stability of DDE, for example, in [7], the method of steps together with the inverse Laplace transform were used to study the stability of a special case of (0.1) when the dimension of $x(t)$ is one.

## 2.1 A simple illustrative example

Consider a simple case of delay differential equation as below

$$\frac{dx(t)}{dt} = ax(t) + bx(t-\delta t) \qquad (2.10)$$

Where in (2.10) $x(t)$ is a 1 dimensional vector, $a, b$ are parameters and $\delta t$ is a small positive number denotes the time delay.

It's easy to see that (2.10) is a special case of (0.1) where $n = N = K = 1$. Now we approximate

$\dfrac{dx(t)}{dt}$ as,

$$\frac{dx(t)}{dt} \approx \frac{x(t)-x(t-\delta t)}{\delta t} \tag{2.11}$$

Then substituting (2.11) into (2.10) yields,

$$\begin{cases} \dfrac{x(t)-x(t-\delta t)}{\delta t} = ax(t)+bx(t-\delta t) \\ \dfrac{x(t-\delta t)-x(t-2\delta t)}{\delta t} = ax(t-\delta t)+bx(t-2\delta t) \end{cases} \tag{2.12}$$

According to the definition of $u(t)$, we have $u(t) = \begin{bmatrix} x(t), x(t-\delta t) \end{bmatrix}^T$ and $u(t-\delta t) = \begin{bmatrix} x(t-\delta t), x(t-2\delta t) \end{bmatrix}^T$. Rewrite (2.12) into a matrix form gives,

$$\frac{1}{\delta t}\begin{bmatrix} 1 & -1 \\ 0 & 1 \end{bmatrix} u(t) + \frac{1}{\delta t}\begin{bmatrix} 0 & 0 \\ -1 & 0 \end{bmatrix} u(t-\delta t) = \begin{bmatrix} a & b \\ 0 & a \end{bmatrix} u(t) + \begin{bmatrix} 0 & 0 \\ 0 & b \end{bmatrix} u(t-\delta t) \tag{2.13}$$

Therefore, in this simple example, we have,

$$D = \begin{bmatrix} 1 & -1 \\ 0 & 1 \end{bmatrix}, \quad D^{(r)} = \begin{bmatrix} 0 & 0 \\ -1 & 0 \end{bmatrix}, \quad B(\theta) = \begin{bmatrix} a & b \\ 0 & a \end{bmatrix} \text{ and } A(\theta) = \begin{bmatrix} 0 & 0 \\ 0 & b \end{bmatrix},$$ this serves as a concrete example to illustrate (2.1).

## 3. Stability analysis of (2.1)

According to the discussion as presented in the previous section, DDE (2.1) is stable if and only if (2.9) holds, to this end, there are 3 cases to be discussed.

i. $DD^T - D^{(r)T}D^{(r)} \succ 0$, in this case, since $\delta t$ can be chosen arbitrarily small, (2.9) holds independent of $\theta$, in other words, (2.1) is stable for all $\theta$.

ii. $DD^T - D^{(r)T}D^{(r)}$ has at least a negative eigenvalue, using similar reasoning as in previous case, we assert that (2.1) is unstable for all $\theta$ in this case.

iii. $DD^T - D^{(r)T}D^{(r)} \succeq 0$, and let matrix $V$ denotes its null space, then (2.1) is stable if and only if below (3.1) holds.

$$V^T E(\theta) V \succ 0 \tag{3.1}$$

Therefore in case iii), the problem of finding $\theta$ such that the DDE system (0.1) stable is equivalent to finding $\theta$, so that (3.1) hold, notice that in (3.1) $V$ is constant matrix, hence verifying the stability of (0.1) can be achieved by solving a semidefinite programming (SDP) like

(3.2) below,

$$\text{minimize } t \text{ s.t. } V^T E(\theta) V \succ tI, t \geq 0 \tag{3.2}$$

where in (3.2), it takes $\theta, t$ as variables.

An important remaining question is how to distinguish the 3 cases of $DD^T - D^{(r)T}D^{(r)}$ as mentioned above? Do they depend on the size of $D$?

The below theorem 3.1 answers this question by showing that the properties of $DD^T - D^{(r)T}D^{(r)}$ depend only on values of $w_1, \ldots, w_m$. As a consequence, one can easily verify that if we use (2.11) as the descrtized version of the derivative operation, then $m = 2$ and $w_1, \ldots, w_m = 1, -1$. Simple calculations reveals that $\tilde{D}^{(t)}\tilde{D}^{(t)T} - \tilde{D}^{(r)T}\tilde{D}^{(r)} = 0$ in theorem 3.1, therefore we have $DD^T - D^{(r)T}D^{(r)} \succeq 0$ in this case.

**Theorem 3.1** If $w_1 \neq 0$, then $DD^T - D^{(r)T}D^{(r)} \succ 0$ if and only if $\tilde{D}^{(t)}\tilde{D}^{(t)T} - \tilde{D}^{(r)T}\tilde{D}^{(r)} \succ 0$, $DD^T - D^{(r)T}D^{(r)} \succeq 0$ if and only if $\tilde{D}^{(t)}\tilde{D}^{(t)T} - \tilde{D}^{(r)T}\tilde{D}^{(r)} \succeq 0$ and $DD^T - D^{(r)T}D^{(r)}$ is indefinite if and only if $\tilde{D}^{(t)}\tilde{D}^{(t)T} - \tilde{D}^{(r)T}\tilde{D}^{(r)}$ is indefinite.

Here in theorem 3.1,

$$\tilde{D}^{(r)} = \begin{bmatrix} w_m & & \\ \vdots & \ddots & \\ w_2 & \cdots & w_m \end{bmatrix}, \text{ with}$$

$$\tilde{D}^{(r)}(j,k) = \begin{cases} w_{m-(k-j)}, & \text{if } k \leq j \\ 0, & \text{otherwise} \end{cases}, 1 \leq j \leq m-1, 1 \leq k \leq m-1 \tag{3.3}$$

And,

$$\tilde{D}^{(t)} = \begin{bmatrix} w_1 & \cdots & w_{m-1} \\ & \ddots & \vdots \\ & & w_1 \end{bmatrix}, \text{ with}$$

$$\tilde{D}^{(t)}(j,k) = \begin{cases} w_{1+(k-j)}, & \text{if } k \geq j \\ 0, & \text{otherwise} \end{cases}, 1 \leq j \leq m-1, 1 \leq k \leq m-1 \tag{3.4}$$

## 4 Extension to higher order differential equations

Consider a higher order differential equations (4.1) as described below,

$$C_1 \frac{d}{dt}x(t) + C_2 \frac{d^2}{dt^2}x(t) + \ldots + C_p \frac{d^p}{dt^p}x(t) = \sum_{i=0}^{K}\sum_{j=0}^{p-1} A_{i,j}(\theta)\frac{d^j x(t-\tau_i)}{dt^j} \quad (4.1)$$

By introducing an extended vector $y(t) := \left[x(t)^T, \frac{d}{dt}x(t)^T, \ldots, \frac{d^p}{dt^p}x(t)\right]^T$, then taking $y(t)$ as variable vector, one obtain a form of 1st order linear delay differentiate equation with respect to $y(t)$, which we state as in (4.2) below:

$$\frac{d}{dt}y(t) = \begin{bmatrix} 0 & I_{(p-1)n} \\ 0 & 0 \end{bmatrix} y(t) + \sum_{i=1}^{K} A_i(\theta) y(t-\tau_i) \quad (4.2)$$

Where in the right hand side of (4.2), $I_{(p-1)n}$ denotes a $(p-1)n \times (p-1)n$ identity matrix, matrices $A_i(\theta), 1 \le i \le K$ are dependent on matrix $A_{i,j}(\theta)$ in (4.1), notices that (4.2) take essentially the form of (0.1), hence one can apply the previous discussed discretizing technique and SDP approach to analyze the stability of (4.2), which is also the stability of (4.1).

## Appendix A numerical differential operator

In this section, we attempt to precisely estimate $f'(t_0)$ with a linear combination of k historical data: $\sum_{i=1}^{k} w_i f(t_0 - i\delta t)$.

According to the Tylor expansion:

$$f(t_0 - i\delta t) = f(t_0) + f'(t_0)(-i\delta t) + \sum_{j=2}^{k} \frac{f^{(n)}(t_0)}{j!}(-i\delta t)^j + o(\delta t)^{k+1}, \quad (A.1)$$
$$i = 1, \ldots, k$$

If the weighting factor $w_i, 1 \le i \le k$ are chosen appropriately such that,

$\sum_{i=1}^{k} w_i = 0$ , $\sum_{i=1}^{k} w_i(-i) \ne 0$ , and $\sum_{i=1}^{k} w_i(-i)^j = 0$ for $j = 2, \ldots, k$ . Then

$\dfrac{\sum_{i=1}^{k} w_i f(t - i\delta t)}{dt \sum_{i=1}^{k} w_i(-i)}$ is an approximation of $f'(t_0)$ with approximation error bounded by

$o(\delta t)^k$. Finally, the full rank property of the following matrix,

$$\begin{bmatrix} 1,\ldots,1 \\ -1,\ldots,-k \\ \cdots \\ (-1)^k,\ldots,(-k)^k \end{bmatrix} \tag{A.2}$$

Implies that such $w_i, 1 \leq i \leq k$ exists, more explicitly, one can choose $\mathrm{w} = [\mathrm{w}_1, \ldots, \mathrm{w}_k]^T$ as,

$$w = c \begin{bmatrix} 1,\ldots,1 \\ -1,\ldots,-k \\ \cdots \\ (-1)^k,\ldots,(-k)^k \end{bmatrix}^{-1} \begin{bmatrix} 1 \\ 0 \\ \vdots \\ 0 \end{bmatrix} \tag{A.3}$$

Where $c \neq 0$ in (A.3) is a constant[2].

# Appendix B proof of theorem 3.1

Notations: For convenience, we adopt the matlab style notations in this section, e.g., if $A$ denotes a matrix, then $A([a:b],[c:d])$ represents one of its sub-matrix constructed by its $a^{th}$, $(a+1)^{th}, \cdots, b^{th}$ rows and $c^{th}, (c+1)^{th}, \cdots, d^{th}$ columns.

According to the structure of $D$ and $D^{(r)}$ as described in (2.2) and (2.3), to prove theorem 3.1, it's sufficient to show $D_i D_i^T - D_i^{(r)T} D_i^{(r)} \succeq 0$ if and only if $\tilde{D}^{(t)} \tilde{D}^{(t)T} - \tilde{D}^{(r)T} \tilde{D}^{(r)} \succeq 0$, see (3.2) and (3.3) for the definitions of matrices $\tilde{D}^{(t)}$ and $\tilde{D}^{(r)}$. For notational convenience, let $\mathcal{G}_N = D_i D_i^T$ denotes a $N \times N$ matrix, with $D_i$ is a matrix defined in (2.2) whose size is also $N \times N$. Similarly, let $\mathcal{H}_N = D_i D_i^T - D_i^{(r)T} D_i^{(r)}$ denote a $N \times N$ matrix, with $D_i$ and $D_i^{(r)}$ defined in (2.2) and (2.3) are matrices of size $N \times N$. Then we have,

---

[2] Following the same technique as introduced in this section, one can also approximate $f'(t_0)$ with a linear combination of $f(t_0)$ and $k-1$ of its historical data $f(t_0 - i\delta t), 1 \leq i \leq k-1$, which reads:

$$f'(t_0) = \frac{\sum_{i=0}^{k-1} w_i f(t - i\delta t)}{dt \sum_{i=0}^{k-1} w_i(-i)} + o(\delta t)^k$$, where coefficients $w_i, 0 \leq i \leq k-1$ are solved similarly as in (A.3).

$$\mathcal{H}_N = \begin{bmatrix} \mathcal{A} & [\mathcal{B},0] \\ [\mathcal{B},0]^T & \mathcal{G}_{N-(m-1)} \end{bmatrix} \quad (B.1)$$

Here in (B.1), $\mathcal{A} = \mathcal{C} - \mathcal{Z}$ is a $(m-1)\times(m-1)$ matrix, with $\mathcal{C} = D_i D_i^T([1:m-1],[1:m-1])$ is a submatrix of $D_i D_i^T$, and $\mathcal{Z} = D_i^{(r)T} D_i^{(r)}([1:m-1],[1:m-1])$ is a submatrix of $D_i^{(r)T} D_i^{(r)}$. $\mathcal{B} = D_i D_i^T([1:m-1],[m:2m-2])$ in (B.1) is a $(m-1)\times(m-1)$ matrix. If we can show that $\mathcal{A}$ is positive definite and the below (B.2) holds,

$$\mathcal{B}\mathcal{A}^{-1}\mathcal{B}^T = \mathcal{Z} \quad (B.2)$$

Then by Schur complement, $\mathcal{H}_N$ is congruent[3] ( see definition 4.5.4 in [8]) to the below matrix,

$$\begin{bmatrix} \mathcal{A} & 0 \\ 0 & \mathcal{H}_{N-(m-1)} \end{bmatrix} \quad (B.3)$$

Without losing generality we assume that $N = k(m-1)$ where $k > 1$ is some positive integer, then apply Schur complement to $\mathcal{H}_{N-(m-1)}$ in (B.3) and iterate this until we hit the bottom $(m-1)\times(m-1)$ sub-matrix of (B.3), we assert that $\mathcal{H}_N$ is congruent to the below matrix,

$$\begin{bmatrix} I & \\ & \mathcal{T}-\mathcal{Z} \end{bmatrix} \quad (B.4)$$

Where in (B.4), $I$ is a $(N-m+1)\times(N-m+1)$ identity matrix. $\mathcal{T} = D_i D_i^T([N-m+2:N],[N-m+2:N])$ is a sub-matrix of $D_i D_i^T$. Therefore at this stage, we assert that $\mathcal{H}_N \succeq 0$ if and only if $\mathcal{T} - \mathcal{Z} \succeq 0$, finally, noticing that $\mathcal{T} = \tilde{D}^{(t)} \tilde{D}^{(t)T}$ and $\mathcal{Z} = \tilde{D}^{(r)T} \tilde{D}^{(r)}$, which proves the conclusion of theorem 3.1.

The remaining part of the proof is to show that (B.2) holds, which we state formally as in below lemma B.1 and then we prove it.

**Lemma B.1** If $w_1 \neq 0$, then $\mathcal{A}$ is positive definite and $\mathcal{B}\mathcal{A}^{-1}\mathcal{B}^T = \mathcal{Z}$ holds.

**Proof:**

---

[3] Let tow $n \times n$ real symmetry matrices A and B be given, we say A and B are congruent (or A is congruent to B) if there exist a non-singular matrix $S$ such that $B = SAS^T$.

To prove the lemma, by Schur complement, it is sufficient to show that below 2 matrices as stated in (B.5) are congruent,

$$\begin{bmatrix} \mathcal{A} & \mathcal{B} \\ \mathcal{B}^T & \mathcal{Z} \end{bmatrix} \sim \begin{bmatrix} I & \\ & 0 \end{bmatrix} \quad \text{(B.5)}$$

Where in the right hand side of (B.5), $I$ is a $(m-1) \times (m-1)$ identity matrix. To show (B.5), firstly we investigate the expressions of elements in matrices $\begin{bmatrix} \mathcal{A} & \mathcal{B} \\ \mathcal{B}^T & 0 \end{bmatrix}$ and $\mathcal{Z}$, for notational simplicity let $\mathcal{U}0$ denotes matrix $\begin{bmatrix} \mathcal{A} & \mathcal{B} \\ \mathcal{B}^T & 0 \end{bmatrix}$, and $\mathcal{U}$ denotes matrix $\begin{bmatrix} \mathcal{A} & \mathcal{B} \\ \mathcal{B}^T & \mathcal{Z} \end{bmatrix}$.

After some elementary calculation, we found,

$$\mathcal{U}0(1, j) = \begin{cases} w(1)w(j), & \text{if } 1 \leq j \leq m \\ 0, & \text{otherwise} \end{cases} \quad \text{(B.6)}$$

$$\mathcal{U}0(i, j) = \begin{cases} \mathcal{U}0(i-1, j-1) + w(i)w(j), & \text{if } 1 < i \leq m-1, 1 < j \leq m \\ \mathcal{U}0(i-1, j-1), & \text{if } 1 < i \leq m-1, m \leq j \leq 2(m-1) \end{cases} \quad \text{(B.7)}$$

$$\mathcal{U}0(i, j) = 0, \text{ if } \min(i, j) \geq m \quad \text{(B.8)}$$

And,

$$Z(i, j) = Z(i+1, j+1) + w(i+1)w(j+1), \text{ if } \max(i, j) < m-1 \quad \text{(B.9)}$$

$$Z(m-1, j) = w(m)w(j+1), 1 \leq j \leq m-1 \quad \text{(B.10)}$$

For notational simplicity, in (B.6)~(B.10) we denote $w = [w_1, \ldots, w_m]^T$ as a m-dimensional vector. As a simple illustrative example, we show the matrix $\mathcal{U}$ in a special case when $m = 3$ as in (B.8) below,

$$\begin{bmatrix} w_1^2, & w_1 w_2, & w_1 w_3, & 0 \\ w_1 w_2, & w_1^2 + w_2^2, & w_1 w_2 + w_2 w_3, & w_1 w_3 \\ w_1 w_3, & w_1 w_2 + w_2 w_3, & w_2^2 + w_3^2, & w_2 w_3 \\ 0, & w_1 w_3, & w_2 w_3, & w_3^2 \end{bmatrix} \quad \text{(B.11)}$$

Then after sequentially applying Schur complement of the 1st diagonal element in $\mathcal{U}$, and then of the 1st diagonal element in $\mathcal{U}(2:2(m-1), 2:2(m-1))$, ..., of the 1st diagonal element in $\mathcal{U}(m-2:2(m-1), m-2:2(m-1))$, we have that $\mathcal{U}$ is congruent to the matrix as stated in (B.12) below,

$$\begin{bmatrix} w_1^2 I & \\ & \mathcal{W} \end{bmatrix} \quad \text{(B.12)}$$

Where in (B.12), $I$ is a $(m-2)\times(m-2)$ identity matrix, $\mathcal{W}=ww^T$ is a rank one $m \times m$ matrix. Finally, applying Schur complement of the 1$^{st}$ diagonal element in $\mathcal{W}$ gives that $\mathcal{U}$ is congruent to $\begin{bmatrix} w_1^2 I & \\ & 0 \end{bmatrix}$, here $I$ is $(m-1)\times(m-1)$ identity matrix. At this stage, we are well prepared to prove the conclusion of the lemma, first, we assert that $\mathcal{A}$ is positive definite, because the preceding Schur complement operations show that $\mathcal{A}$ is congruent to $w_1^2 I$. Secondly, we assert that $\mathcal{B}\mathcal{A}^{-1}\mathcal{B}^T = \mathcal{Z}$, because otherwise if $\mathcal{B}\mathcal{A}^{-1}\mathcal{B}^T \neq \mathcal{Z}$, after applying Schur complement of $\mathcal{A}$ in matrix $\mathcal{U} = \begin{bmatrix} \mathcal{A} & \mathcal{B} \\ \mathcal{B}^T & 0 \end{bmatrix}$, we have that $\mathcal{U}$ is congruent to $\begin{bmatrix} \mathcal{A} & \\ & \mathcal{Z}-\mathcal{B}\mathcal{A}^{-1}\mathcal{B}^T \end{bmatrix}$, but according to the preceding discussion, $\mathcal{U}$ is also congruent to $\begin{bmatrix} w_1^2 I & \\ & 0 \end{bmatrix}$ and $\mathcal{A}$ is congruent to $w_1^2 I$, the fact that $\mathcal{B}\mathcal{A}^{-1}\mathcal{B}^T \neq \mathcal{Z}$ would imply that $\mathcal{U}$ is congruent to 2 matrices with different inertias[4] (definition 4.5.6 [8]), which is clearly a contradiction and therefore the conclusion of lemma B.1 is proved.∎

---

[4] Let a $n \times n$ matrix A be given, then the inertia of A is defined as the order triple: $i(A)=(i_+(A), i_-(A), i_0(A))$, where $i_+(A)$ is the number of positive eigenvalues of A, $i_-(A)$ is the number of negative eigenvalues of A and $i_0(A)$ is the number of zero eigenvalues of A.